\documentclass[11pt]{article}
\usepackage{graphicx,subfigure}
\usepackage{epsfig}
\usepackage{amssymb}
\usepackage{mathrsfs}
\usepackage{color}
\usepackage[dvipdfm]{hyperref}

\newcommand{\E}{{\cal E}}

\newtheorem{theorem}{Theorem}[section]
\newtheorem{definition}{Definition}[section]
\newtheorem{lemma}[theorem]{Lemma}

\newtheorem{example}[theorem]{Example}

\def\whitebox{{\hbox{\hskip 1pt
 \vrule height 6pt depth 1.5pt
 \lower 1.5pt\vbox to 7.5pt{\hrule width
    3.2pt\vfill\hrule width 3.2pt}%
 \vrule height 6pt depth 1.5pt
 \hskip 1pt } }}
\def\qed{\ifhmode\allowbreak\else\nobreak\fi\hfill\quad\nobreak
     \whitebox\medbreak}

\newcommand{\ignore}[1]{}

\begin {document}

\baselineskip 16pt
\title{Resistance distances in corona and neighborhood corona graphs with Laplacian generalized inverse approach}

 \author{\small   Jia\textrm{-}Bao \ Liu$^{a,b}$, \ \  Xiang\textrm{-}Feng \ Pan$^{a,}
 $\thanks{Corresponding author. Tel:+86-551-63861313. \E-mail:liujiabaoad@163.com (J.Liu),
  xfpan@ahu.edu.cn(X.Pan), hufu@mail.ustc.edu.cn(F.Hu).} ,\ \ Fu\textrm{-}Tao Hu$^{a}$\\
 \small  $^{a}$ School of Mathematical Sciences, Anhui  University, Hefei 230601, P. R. China\\
 \small  $^{b}$ Department of Public Courses, Anhui Xinhua
 University, Hefei 230088, P. R. China\\}

\date{}
\maketitle
\begin{abstract}
  Let $G_1$ and $G_2$ be two graphs on disjoint sets of $n_1$ and
$n_2$ vertices, respectively. The corona of graphs $G_1$ and
$G_2$, denoted by $G_1\circ G_2$, is the graph formed from one
copy of $G_1$ and $n_1$ copies of $G_2$ where the $i$-th vertex of
$G_1$ is adjacent to every vertex in the $i$-th copy of $G_2$. The
neighborhood corona of $G_1$ and $G_2$, denoted by $G_1\diamond
G_2$, is the graph obtained by taking one copy of $G_1$ and $n_1$
copies of $G_2$ and joining every neighbor of the $i$-th vertex of
$G_1$ to every vertex in the $i$-th copy of $G_2$ by a new edge.
In this paper, the Laplacian generalized inverse for the
 graphs $G_1\circ G_2$ and $G_1\diamond
G_2$ are investigated, based on which the resistance distances of
any two vertices in $G_1\circ G_2$ and $G_1\diamond G_2$ can be
obtained. Moreover, some examples as applications are presented,
which illustrate the correction and efficiency of the proposed
method.

\end{abstract}

\section{ Introduction}

All graphs considered in this paper are simple and undirected. Let
$G = (V(G),E(G))$ be a graph with vertex set $V(E) = \{v_1, v_2,
\dots , v_n\}$ and edge set $E(G) = \{e_1, e_2, \dots , e_m\}$.
The adjacency matrix of $G$, denoted by $A(G)$, is the $n \times
n$ matrix whose $(i,j)$-entry is $1$ if $v_i$ and $v_j$ are
adjacent in $G$ and $0$ otherwise. Denote $D(G)$ to be the
diagonal matrix with diagonal entries $d_G(v_1),d_G(v_2),\dots,
d_G(v_n).$ The Laplacian matrix of $G$ is defined as $L(G)=D(G)-
A(G)$. For other undefined notations and terminology from graph
theory, the readers may refer to~\cite{Bondy1976} and the
references therein.

 The conventional distance between vertices $v_i$ and $v_j$, denoted by
$d_{i, j}$, is the length of a shortest path between them.
 Klein and Randi\'c~\cite{Klein1993} introduced a new distance function named
resistance distance based on electrical network theory, the
resistance distance between vertices $i$ and $j$, denoted by
$r_{ij}$, is defined to be the effective electrical resistance
between them if each edge of $G$ is replaced by a unit
resistor~\cite{Klein1993}.
   For more information on resistance distance
of graphs, the readers are referred to the most recent papers
~\cite{Bu2014,Y2013,Feng2014,Ya2014,YangJ2014,Zhang2013,Liu2014,Liu2015,LiuP2015}.

Until now, many graph operations such as the Cartesian product,
the Kronecker product, the corona and neighborhood corona graphs
have been introduced
in~\cite{Liu2013,Lu2013,McLeman2011,Wang2012,Gopalapillai2011}.
Let $G_1$ and $G_2$ be two vertex disjoint graphs. The following
definition comes from~\cite{Lu2013}.

\begin{definition}\label{1-1}(see~\cite{Lu2013})
Let $G_1$ and $G_2$ be two graphs on disjoint sets of $n_1$ and
$n_2$ vertices, respectively. The corona of two graphs $G_1$ and
$G_2$ is the graph $G = G_1\circ G_2$ formed from one copy of
$G_1$ and $n_1$ copies of $G_2$ where the $i$-th vertex of $G_1$
is adjacent to every vertex in the $i$-th copy of $G_2$.
\end{definition}

The neighborhood corona, which is a variant of the corona
operation, was recently introduced in~\cite{Gopalapillai2011}.

\begin{definition}\label{1-2}(see~\cite{Gopalapillai2011})
Let $G_1$ and $G_2$ be two graphs on disjoint sets of $n_1$ and
$n_2$ vertices, respectively. The neighborhood corona of $G_1$ and
$G_2$, denoted by $G_1\diamond G_2$, is the graph  obtained by
taking one copy of $G_1$ and $n_1$ copies of $G_2$ and joining
every neighbor of the $i$-th vertex of $G_1$ to every vertex in
the $i$-th copy of $G_2$ by a new edge.
\end{definition}

\begin{figure}[ht]
\center
  \includegraphics[width=\textwidth]{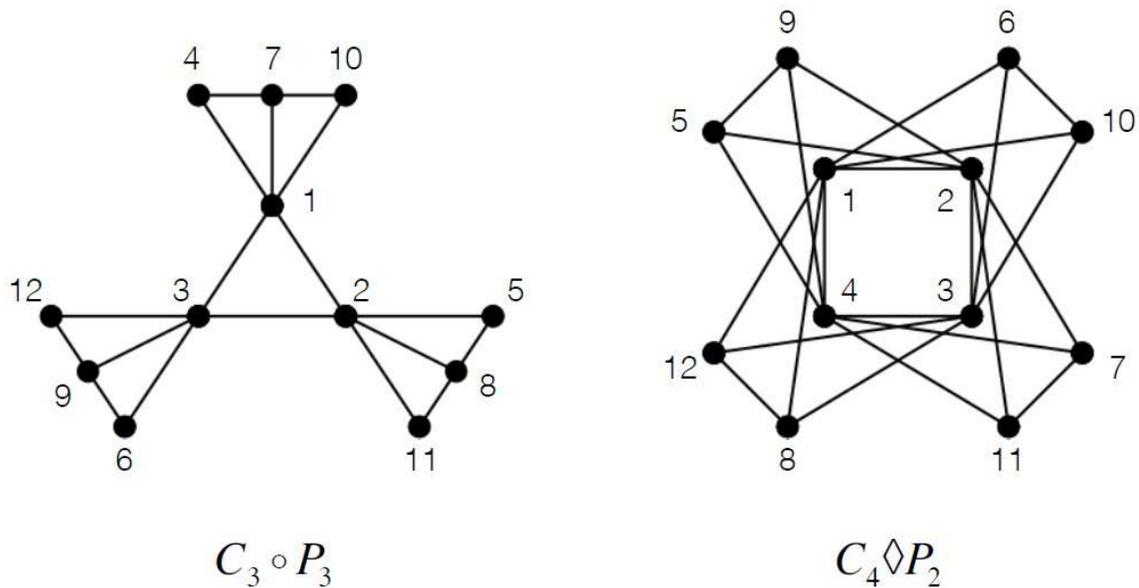}
  \caption{ $C_3\circ P_3$ and $C_4\diamond P_2$.}
  \vspace{+1em}
\end{figure}

Let $P_n$ and  $C_n$ denote a path and cycle with $n$ vertices,
respectively. From the definitions,  Figure 1 shows the graphs
$C_3\circ P_3$ and $C_4\diamond P_2$.

 Bu et al. investigated
 resistance distances in subdivision-vertex join and subdivision-edge join of graphs~\cite{Bu2014}.
 Motivated by the
results, in this paper, we further explored the Laplacian
generalized inverse for the corona and neighborhood corona graphs,
based on which all the resistance distances between arbitrary two
vertices can be directly obtained via simple calculations.

\section{ Preliminaries and Lemmas }

At the beginning of this section, we review some concepts in
matrix theory.
 Let $A$ be a matrix, $X$ is called the
$\{1 \}$ inverse of $A$ and denoted by $A^{\{1 \}}$, if $X$
satisfies the following condition: $AXA=A.$ Given a square matrix
$A$, the group inverse of $A$, denoted by $A^{\#}$, is the unique
matrix $X $ that satisfies matrix equations~\cite{Bu2014} $ (I).
AXA=A,~(II). XAX=X,~(III). AX=XA.$
 If $A$ is real symmetric, then $A^{\#}$ exists and $A^{\#}$ is a
symmetric $\{1 \}$-inverse of $A$. In fact, $A^{\#}$ is equal to
the Moore-Penrose inverse of $A$ since $A$ is symmetric
~\cite{Bu2014}.

The Kronecker product $A\bigotimes B$ of two matrices $A =
(a_{ij})_{m\times n}$ and $B = (b_{ij})_{p\times q}$ is the
$mp\times nq $ matrix obtained from $A$ by replacing each element
$a_{ij} $ by $a_{ij}B $.  The reader is referred
to~\cite{Horn1991} for other properties of the Kronecker product
not mentioned here.

 It is known that resistance distances in a connected graph $G$
can be obtained from any $\{1\}$-inverse of $L(G)$ according to
the following lemma (see \cite{Bu2014}).

\begin{lemma}\label{2-1}   (see~\cite{Bu2014}) Let $G$ be a connected graph, and $(L_G)_{ij}$ denote the ${(i,j)}$-entry of $L_G$.
Then
$$ r_{ij}(G) =(L_G^{(1)})_{ii} + (L_G^{(1)})_{jj}- (L_G^{(1)})_{ij}- (L_G^{(1)})_{ji}
= (L_G^{\#})_{ii} + (L_G^{\#})_{jj}- 2(L_G^{\#})_{ij} .$$
\end{lemma}

\begin{lemma}\label{2-2} (see \cite{Bu2014}) Let $M=\left[%
\begin{array}{cc}
  A & B \\
  C & D \\
\end{array}%
\right]$ be a nonsingular matrix. If $A$ and $D$ are nonsingular,
then
$$ M^{-1} =\left[%
\begin{array}{cc}
  A^{-1}+A^{-1}BS^{-1}CA^{-1} & ~~~~-A^{-1}BS^{-1} \\
  -S^{-1}CA^{-1} &~~~~ S^{-1} \\
\end{array}%
\right],\\ $$

 where $ S=D-CA^{-1}B$ is the Schur complement of $A$ in $M$.
\end{lemma}

The following similar result holds for Laplacian matrix of a
connected graph.

\begin{lemma}\label{2-3} (see \cite{Zhou2014})   Let $L=\left[%
\begin{array}{cc}
  L_1 &  L_2 \\
  L_2^T &  L_3 \\
\end{array}%
\right]$ be the Laplacian matrix of a connected graph. If $L_1$ is
nonsingular, then
$X=\left[%
\begin{array}{cc}
  L_1^{-1}+L_1^{-1}L_2S^{\#}L_2^TL_1^{-1} &~~~ -L_1^{-1}L_2S^{\#} \\
  -S^{\#}L_2^TL_1^{-1} &~~~ S^{\#} \\
\end{array}%
\right]$ is a symmetric $\{1\}$-inverse of $L$, where
$S=L_3-L_2^TL_1^{-1}L_2.$
\end{lemma}

\section{\bf The Laplacian generalized inverse for
 graphs $G_1\circ G_2$ and $ G_1\diamond G_2$}

\subsection{\bf The Laplacian generalized inverse for
 graph $G_1\circ G_2$}

Let $\mathbf{1}_n$ and $J_{n\times n}$ be all-one column vector of
dimensions $n$ and all-one $n \times n$ matrix, respectively.

 {\bf
Theorem 3.1} Let $G_1$ be an $r_1$-regular graph with $n_1$
vertices and $m_1$ edges, and $G_2$ an arbitrary graph with $n_2$
vertices, then the following matrix
$$\left[%
\begin{array}{c|lr}
   L_1^{-1}+L_1^{-1}L_2S^{\#}L_2^TL_1^{-1}  &~~ -L_1^{-1}L_2S^{\#}\\
   \hline
  -S^{\#}L_2^TL_1^{-1}     & ~~~~~~~S^{\#} \\
\end{array}%
\right]$$

is a symmetric $\{1\}$-inverse of $L_{(G_1\circ G_2)}$, where

$L_1=[L_{(G_1)}+n_2I_{n_1}],
 L_2=[-\mathbf{1}^T_{n_2}\otimes I_{n_1}],
 L_3=(L(G_2)+I_{n_2})\otimes I_{n_1},
 S=\Big[L_3 -J_{n_2\times n_2} \otimes L_1^{-1}
\Big].$

{\bf Proof.} Let $G_1$ be an arbitrary $r_1$-regular graphs with
$n_1$ vertices and $m_1$ edges, and $G_2$ an arbitrary graphs with
$n_2$ vertices, respectively. Label the vertices of $G_1\circ G_2$
as follows. Let  $V(G_1) = \{v_1, v_2,\dots , v_{n_1}\}$  and
$V(G_2) = \{w_1, w_2,\dots , w_{n_2}\}$. For $i = 1, 2, \dots ,
n_1$, let $w^i_1, w^i_2, \dots , w_{n_2}^i$ denote the vertices of
the $i$-th copy of $G_2$, with the understanding that $w^i_j$ is
the copy of $w_j$ for each $j$. Denote $W_j =\{ w^1_j, w^2_j,\dots
, w^{n_1}_j\}$, for $j = 1, 2,\dots , n_2.$ Then
\begin{equation}\label{}
V(G_1)  \bigcup \big[W_1 \bigcup W_2 \bigcup \dots  \bigcup
W_{n_2} \big]
\end{equation}
 is a partition
of $V (G_1\circ G_2)$. Obviously, the degrees of the vertices of
$G_1\circ G_2$ are:
 $d_{G_1\circ G_2}(e_i) = 2 $, for $i = 1, 2,
\dots , m_1$, $d_{G_1\circ G_2}(v_i) = n_2 + d_{G_1}(v_i)$, for $i
= 1, 2, \dots , n_1$, and $d_{G_1\circ G_2}(w_j^i) = d_{G_2}(w_j)
+ 1$, for $i = 1, 2, \dots , n_1, j = 1, 2,\dots , n_2$.

 Since $G_1$ is an
$r_1$-regular graph, we have $D(G_1) = r_1I_{n_1} $. With respect
to the partition (1), then the Laplacian matrix of $G_1\circ G_2$
can be written as
$$ L (G_1\circ G_2)=  \left[%
\begin{array}{c|lr}
 L(G_1)+n_2I_{n_1} &~~~ -\mathbf{1}^T_{n_2}\otimes I_{n_1} \\
 \hline
 -\mathbf{1}_{n_2}\otimes I_{n_1}  & ~~(L(G_2)+I_{n_2})\otimes I_{n_1} \\
\end{array}%
\right].$$

We begin with the calculation  $S$. For convenience, let

$L_1=[L_{(G_1)}+n_2I_{n_1}], L_2=[-\mathbf{1}^T_{n_2}\otimes
I_{n_1}],$ $ L_2^T=[-\mathbf{1}_{n_2}\otimes I_{n_1}],
L_3=(L(G_2)+I_{n_2})\otimes I_{n_1}.$\\

By Lemma 2.3, we have
\begin {eqnarray*}
  S &&=\Big[L(G_2)+I_{n_2}\Big]\otimes I_{n_1}-\left[%
\begin{array}{c}
 -\mathbf{1}_{n_2}\otimes I_{n_1}\\
\end{array}%
\right]\left[%
\begin{array}{c}
   L(G_1)+n_2I_{n_1}\\
\end{array}%
\right]^{-1}         \left[%
\begin{array}{c}
 -\mathbf{1}^T_{n_2}\otimes I_{n_1}  \end{array}%
\right]  \\
&&=\Big[L(G_2)+I_{n_2}\Big]\otimes I_{n_1}
-J_{n_2\times n_2} \otimes \left[%
\begin{array}{c}
 L(G_1)+n_2I_{n_1} \\
\end{array}%
\right]^{-1}  \\
&&=L_3 -J_{n_2\times n_2} \otimes L_1^{-1}.
\end {eqnarray*}

 Based on  Lemma 2.3, the following matrix
$$\left[%
\begin{array}{c|lr}
   L_1^{-1}+L_1^{-1}L_2S^{\#}L_2^TL_1^{-1}  & ~-L_1^{-1}L_2S^{\#}\\
   \hline
  -S^{\#}L_2^TL_1^{-1}     & ~~~~~~~~S^{\#} \\
\end{array}%
\right]$$
 is a symmetric $\{1\}$-inverse of $L_{(G_1\circ G_2)}$,
where

$L_1=[L_{(G_1)}+n_2I_{n_1}],
 L_2=[-\mathbf{1}^T_{n_2}\otimes I_{n_1}],
 L_3=(L(G_2)+I_{n_2})\otimes I_{n_1},
 S=\Big[L_3 -J_{n_2\times n_2} \otimes L_1^{-1}
\Big].$ \hfill$\blacksquare$

\subsection{\bf The Laplacian generalized inverse for
 graph $G_1\diamond G_2$}

When $G_1$ is a regular graph, we obtain the Laplacian generalized
inverse for graph
 $ G_1\diamond G_2$ as follows.

{\bf Theorem 3.2} Let $G_1$ be an $r_1$-regular graph with $n_1$
vertices and $m_1$ edges, and $G_2$ an arbitrary graph with $n_2$
vertices, then the following matrix

$$\left[%
\begin{array}{c|lr}
   L_1^{-1}+L_1^{-1}L_2S^{\#}L_2^TL_1^{-1}  & ~-L_1^{-1}L_2S^{\#}\\
   \hline
  -S^{\#}L_2^TL_1^{-1}     & ~~~~~~~~S^{\#} \\
\end{array}%
\right]$$
 is a symmetric $\{1\}$-inverse of $L_{(G_1\diamond
G_2)}$, where
$L_1= L(G_1)+ n_2D(G_1), L_2=-\mathbf{1}_{n_2}^T\otimes A(G_1),$\\
$ L_3=L(G_2)\otimes I_{n_1}+I_{n_2} \otimes
D(G_1),S=L_3-J_{n_2\times n_2} \otimes \big[ A(G_1)^T L_1^{-1}
A(G_1) \big].$

{\bf Proof.} We label the vertices of $ G_1\diamond G_2$ as
follows. Let  $V(G_1) = \{v_1, v_2,\dots ,v_{n_1}\}$, and $V(G_2)
= \{w_1, w_2,\dots , w_{n_2}\}$. For $i = 1, 2, \dots , n_1$, let
$w^i_1, w^i_2, \dots , w_{n_2}^i$ denote the vertices of the
$i$-th copy of $G_2$, with the understanding that $w^i_j$ is the
copy of $w_j$ for each $j$. Denote $U_j =\{ w^1_j, w^2_j,\dots ,
w^{n_1}_j\}$, for $j = 1, 2,\dots , n_2.$ Then
\begin{equation}\label{}
V(G_1)\bigcup \Big [U_1 \bigcup U_2 \bigcup \dots \bigcup
U_{n_2}\Big]
\end{equation}
 is a partition of $V
(G_1\diamond G_2)$. Clearly, the degrees of the vertices of
$G_1\diamond G_2$ are:

 $d_{G_1\diamond G_2}(v_i) =
(n_2+1)d_{G_1}(v_i),$ ~for $i = 1, 2, \dots ,n_1$, and
$d_{G_1\diamond G_2}(w_j^i) = d_{G_2}(w_j) + d_{G_1}(v_i),$ ~for
$i = 1, 2, \dots , m_1, j = 1, 2,\dots , n_2$.

 Since $G_1$ is an $r_1$-regular graph, we have $D(G_1) =
r_1I_{n_1} $.  With respect to the partition (2), then the
Laplacian matrix of $ G_1\diamond G_2$ can be written as

$$ L (G_1\diamond G_2)= \left[%
\begin{array}{c|lr}
~L(G_1)+ n_2D(G_1)~ & ~~~~~~~~~-\mathbf{1}_{n_2}^T\otimes A(G_1) \\
  \hline
~-\mathbf{1}_{n_2}\otimes A(G_1)^T  & ~~L(G_2)\otimes I_{n_1}+I_{n_2} \otimes D(G_1) \\
\end{array}%
\right].$$

For convenience, let $L_1= L(G_1)+ n_2D(G_1),
L_2=-\mathbf{1}_{n_2}^T\otimes A(G_1),
L_2^T=-\mathbf{1}_{n_2}\otimes A(G_1)^T,$\\ $
L_3=L(G_2)\otimes I_{n_1}+I_{n_2} \otimes D(G_1).$\\

Similarly, by Lemma 2.3, we have
\begin {eqnarray*}
  S &&=\big[L(G_2)\otimes I_{n_1}+I_{n_2} \otimes D(G_1)\big] -\left[%
-\mathbf{1}_{n_2}\otimes A(G_1)^T\right]
\left[%
\begin{array}{c}
  L(G_1)+ n_2D(G_1)
\end{array}%
\right]^{-1}         \left[%
\begin{array}{c}
 -\mathbf{1}_{n_2}^T\otimes A(G_1)  \end{array}%
\right]  \\
&&= L_3-J_{n_2\times n_2} \otimes
\big[ A(G_1)^T L_1^{-1} A(G_1) \big]  . \\
\end {eqnarray*}
 Based on  Lemma 2.3, the following matrix
$$\left[%
\begin{array}{c|lr}
   L_1^{-1}+L_1^{-1}L_2S^{\#}L_2^TL_1^{-1}  & ~-L_1^{-1}L_2S^{\#}\\
   \hline
  -S^{\#}L_2^TL_1^{-1}     & ~~~~~~~~S^{\#} \\
\end{array}%
\right]$$ is a symmetric $\{1\}$-inverse of $L_{(G_1\circ G_2)}$,
where $L_1= L(G_1)+ n_2D(G_1), L_2=-\mathbf{1}_{n_2}^T\otimes
A(G_1),$\\ $L_3=L(G_2)\otimes I_{n_1}+I_{n_2} \otimes D(G_1),
S=L_3-J_{n_2\times n_2} \otimes \big[ A(G_1)^T L_1^{-1} A(G_1)
\big].$
 \hfill$\blacksquare$

\section{Applications and some examples}

 As an application of the proposed theorems, we present some examples to show all the resistance distances of
 any two vertices in graphs $G_1\circ G_2$ and $G_1\diamond G_2$ can be obtained by the proposed method.

\begin{example}\label{4-1}
Laplacian generalized inverse for $ C_3\circ P_3$ and resistance
distances matrix.
\end{example}\label{4-1}

The Laplacian matrix
$ L_{(C_3\circ P_3)}=  \left[%
\begin{array}{c|lr}
 L(C_3)+3I_{3} &~~~~ -\mathbf{1}^T_{3}\otimes I_{3} \\
 \hline
 -\mathbf{1}_{3}\otimes I_{3}  & ~(L(P_3)+I_{n_2})\otimes I_{3} \\
\end{array}%
\right].$

 Based on Theorem 3.1, we can obtain that

$$L^{\{1\}}_{(C_3\circ P_3)}=
\left[%
\begin{array}{cccccccccccc}
\frac{1}{3} & 0 & 0 & \frac{2}{9} & \frac{-1}{9}  & \frac{-1}{9} &\frac{2}{9} & \frac{-1}{9} & \frac{-1}{9} & \frac{2}{9} &\frac{-1}{9} &\frac{-1}{9} \\
0 &\frac{1}{3} & 0 & \frac{-1}{9} & \frac{2}{9} &\frac{-1}{9}& \frac{-1}{9} &\frac{2}{9} & \frac{-1}{9} & \frac{-1}{9}& \frac{2}{9} & \frac{-1}{9}\\
0 & 0  & \frac{1}{3} &\frac{-1}{9}  & \frac{-1}{9} &\frac{2}{9} &\frac{-1}{9}  & \frac{-1}{9} &\frac{2}{9} & \frac{-1}{9}  & \frac{-1}{9} &\frac{2}{9}\\
\frac{2}{9} & \frac{-1}{9}  & \frac{-1}{9} & \frac{53}{72} &\frac{-2}{9} & \frac{-2}{9} & \frac{ 13}{36} & \frac{-2}{9} &\frac{-2}{9} &    \frac{ 17}{72} & \frac{-2}{9} & \frac{-2}{9}\\
 \frac{-1}{9}  &\frac{2}{9} & \frac{-1}{9} & \frac{-2}{9} & \frac{53}{72} &\frac{-2}{9}  & \frac{-2}{9} & \frac{ 13}{36} &\frac{-2}{9} & \frac{-2}{9} & \frac{ 17}{72} &\frac{-2}{9}\\
 \frac{-1}{9}  & \frac{-1}{9} & \frac{2}{9} &  \frac{-2}{9} & \frac{-2}{9} &  \frac{53}{72} &   \frac{-2}{9}  & \frac{-2}{9} & \frac{ 13}{36} & \frac{-2}{9} &\frac{-2}{9} &    \frac{ 17}{72}\\
\frac{2}{9} & \frac{-1}{9}  & \frac{-1}{9} &  \frac{ 13}{36}
&\frac{-2}{9} &\frac{-2}{9} & \frac{11}{18}    & \frac{-2}{9}
&\frac{-2}{9} &  \frac{ 13}{36} &
\frac{-2}{9} &\frac{-2}{9}\\
\frac{-1}{9}  &\frac{2}{9} & \frac{-1}{9}  &   \frac{-2}{9} &\frac{ 13}{36}  &   \frac{-2}{9}  & \frac{-2}{9} & \frac{11}{18}    & \frac{-2}{9}  & \frac{-2}{9} & \frac{ 13}{36} & \frac{-2}{9}\\
\frac{-1}{9}  & \frac{-1}{9} & \frac{2}{9} &   \frac{-2}{9}&\frac{-2}{9} &  \frac{ 13}{36}  &  \frac{-2}{9}  & \frac{-2}{9} &\frac{11}{18}  &  \frac{-2}{9}  & \frac{-2}{9}  &\frac{ 13}{36} \\
\frac{2}{9} & \frac{-1}{9}  & \frac{-1}{9} &   \frac{ 17}{72} &\frac{-2}{9} & \frac{-2}{9} & \frac{ 13}{36}  &   \frac{-2}{9}  &\frac{-2}{9} &   \frac{53}{72} &   \frac{-2}{9}  & \frac{-2}{9} \\
\frac{-1}{9}  &\frac{2}{9} & \frac{-1}{9}  & \frac{-2}{9} &
 \frac{ 17}{72} &\frac{-2}{9} & \frac{-2}{9} &  \frac{ 13}{36}  &  \frac{-2}{9} & \frac{-2}{9} &   \frac{53}{72} &   \frac{-2}{9}\\
 \frac{-1}{9}  & \frac{-1}{9} & \frac{2}{9}  &  \frac{-2}{9}  & \frac{-2}{9}  &   \frac{ 17}{72} & \frac{-2}{9}  & \frac{-2}{9}    &\frac{ 13}{36}  & \frac{-2}{9}  & \frac{-2}{9}  &
  \frac{53}{72}\\
\end{array}%
\right].$$

By Lemma 2.1 and $L^{\{1\}}_{(C_3\circ P_3)}$, the resistance
distances matrix of $C_3\circ P_3$ is
$$R_{(C_3\circ P_3)}=
\left[%
\begin{array}{cccccccccccc}
 0 &\frac{2}{3} & \frac{2}{3} & \frac{5}{8} & \frac{31}{24}  & \frac{31}{24} &\frac{1}{2} & \frac{7}{6} & \frac{7}{6} & \frac{5}{8} &\frac{31}{24} &\frac{31}{24}\\
\frac{2}{3} &0 & \frac{2}{3} & \frac{31}{24} & \frac{5}{8} &\frac{31}{24}& \frac{7}{6} &\frac{1}{2} & \frac{7}{6} & \frac{31}{24}& \frac{5}{8} & \frac{31}{24}\\
\frac{2}{3} & \frac{2}{3} & 0 &\frac{31}{24}  & \frac{31}{24} &\frac{5}{8} &\frac{7}{6}  & \frac{7}{6} &\frac{1}{2} & \frac{31}{24}  & \frac{31}{24} &\frac{5}{8}\\
\frac{5}{8} & \frac{31}{24}  & \frac{31}{24} & 0 &\frac{23}{12} & \frac{23}{12} & \frac{5}{8} & \frac{43}{24} &\frac{43}{24} & 1 & \frac{23}{12} & \frac{23}{12}\\
\frac{31}{24}  &\frac{5}{8} & \frac{31}{24} &\frac{23}{12} & 0 &\frac{23}{12}  & \frac{43}{24} & \frac{5}{8} &\frac{43}{24} & \frac{23}{12} & 1 &\frac{23}{12}\\
 \frac{31}{24}  & \frac{31}{24} &\frac{5}{8} &\frac{23}{12} & \frac{23}{12} &0 & \frac{43}{24}  &\frac{43}{24} & \frac{5}{8} & \frac{23}{12} &\frac{23}{12} &  1\\
\frac{1}{2}&\frac{7}{6}&\frac{7}{6}&\frac{5}{8}&\frac{43}{24}&\frac{43}{24}&0&\frac{5}{3}&\frac{5}{3}&\frac{5}{8}&\frac{43}{24} &\frac{43}{24}\\
\frac{7}{6} &\frac{1}{2}&\frac{7}{6}  &\frac{43}{24} &\frac{5}{8}&\frac{43}{24}  &\frac{5}{3} & 0    &\frac{5}{3}  &\frac{43}{24} &\frac{5}{8} &\frac{43}{24}\\
\frac{7}{6}&\frac{7}{6} & \frac{1}{2} &\frac{43}{24}&\frac{43}{24} &  \frac{5}{8}  &\frac{5}{3}  &\frac{5}{3} &0  & \frac{43}{24}  &\frac{43}{24}  &\frac{5}{8}\\
\frac{5}{8} & \frac{31}{24} &\frac{31}{24} &1 &\frac{23}{12} & \frac{23}{12} &\frac{5}{8} &\frac{43}{24}  &\frac{43}{24}  &0& \frac{23}{12} &\frac{23}{12} \\
\frac{31}{24}& \frac{5}{8} & \frac{31}{24}  &\frac{23}{12} &1 &\frac{23}{12} & \frac{43}{24} &\frac{5}{8} &\frac{43}{24} &\frac{23}{12} &0 & \frac{23}{12}\\
\frac{31}{24}  & \frac{31}{24} &\frac{5}{8} & \frac{23}{12} & \frac{23}{12}  & 1 & \frac{43}{24}&\frac{43}{24} &  \frac{5}{8}& \frac{23}{12} & \frac{23}{12} &0\\
\end{array}%
\right],$$
 where $r_{ij}$ denotes resistance distance of two
vertices between $i$ and $j.$\hfill$\blacksquare$

\begin{example}\label{4-2}
Laplacian generalized inverse for $ C_4\diamond P_2$ and
resistance distances matrix.
\end{example}\label{4-2}

Completely similar deduction by Theorem 3.2, we can obtain
$$L^{\{1\}}_{(C_4\diamond P_2)}=
\left[%
\begin{array}{cccccccccccc}
\frac{5}{24} & 0 & \frac{1}{24} & 0 & \frac{-1}{16}  & \frac{1}{16} &\frac{-1}{16}  & \frac{1}{16} & \frac{-1}{16}  & \frac{1}{16} &\frac{-1}{16}  & \frac{1}{16} \\
0 &\frac{5}{24} & 0 & \frac{1}{24} & \frac{1}{16} &\frac{-1}{16}  & \frac{1}{16} & \frac{-1}{16}  & \frac{1}{16} &\frac{-1}{16}  & \frac{1}{16} & \frac{-1}{16} \\
\frac{1}{24} & 0 & \frac{5}{24} & 0 & \frac{-1}{16}  & \frac{1}{16} &\frac{-1}{16}  & \frac{1}{16} & \frac{-1}{16}  & \frac{1}{16} &\frac{-1}{16}  & \frac{1}{16}\\
0 &\frac{1}{24} & 0 & \frac{5}{24} & \frac{1}{16} &\frac{-1}{16}  & \frac{1}{16} & \frac{-1}{16}  & \frac{1}{16} &\frac{-1}{16}  & \frac{1}{16} & \frac{-1}{16}\\
 \frac{-1}{16}  & \frac{1}{16} &\frac{-1}{16}  & \frac{1}{16} & \frac{3}{8} &\frac{-1}{8}  & 0 & \frac{-1}{8} &\frac{1}{8} & \frac{-1}{8} & 0 &\frac{-1}{8}\\
\frac{1}{16} & \frac{-1}{16}  & \frac{1}{16} &\frac{-1}{16} & \frac{-1}{8} &  \frac{3}{8} & \frac{-1}{8}  & 0 & \frac{-1}{8} &\frac{1}{8} & \frac{-1}{8} & 0\\
 \frac{-1}{16}  & \frac{1}{16} &\frac{-1}{16}  & \frac{1}{16}
&0 & \frac{-1}{8} &\frac{3}{8} & \frac{-1}{8} &0 & \frac{-1}{8} &\frac{1}{8} & \frac{-1}{8}\\
\frac{1}{16} & \frac{-1}{16}  & \frac{1}{16} &\frac{-1}{16} &\frac{-1}{8} &0 & \frac{-1}{8} &\frac{3}{8}  & \frac{-1}{8} &0 & \frac{-1}{8} &\frac{1}{8}\\
\frac{-1}{16}  & \frac{1}{16} &\frac{-1}{16}  & \frac{1}{16} &\frac{1}{8} &  \frac{-1}{8} &0 & \frac{-1}{8} &\frac{3}{8}  & \frac{-1}{8}  & 0 & \frac{-1}{8} \\
\frac{1}{16} & \frac{-1}{16}  & \frac{1}{16} &\frac{-1}{16} &\frac{-1}{8} &\frac{1}{8} & \frac{-1}{8} & 0  &\frac{-1}{8} &\frac{3}{8} & \frac{-1}{8} & 0 \\
\frac{-1}{16}  & \frac{1}{16} &\frac{-1}{16}  & \frac{1}{16}&
 0 & \frac{-1}{8} &\frac{1}{8}  & \frac{-1}{8}  &  0 & \frac{-1}{8} &\frac{3}{8}  & \frac{-1}{8}\\
\frac{1}{16} & \frac{-1}{16}  & \frac{1}{16} &\frac{-1}{16} &
\frac{-1}{8} &0 & \frac{-1}{8} &\frac{1}{8} &\frac{-1}{8} &0 & \frac{-1}{8} &\frac{3}{8}\\
\end{array}%
\right].$$

By Lemma 2.1 and $L^{\{1\}}_{(C_4\diamond P_2)}$, the resistance
distances matrix of $C_4\diamond P_2$ is
$$R_{(C_4\diamond P_2)}=
\left[%
\begin{array}{cccccccccccc}
0 & \frac{5}{12} &\frac{1}{3} &  \frac{5}{12} & \frac{17}{24}  & \frac{11}{24} &\frac{17}{24}  & \frac{11}{24} & \frac{17}{24}  & \frac{11}{24} &\frac{17}{24}  & \frac{11}{24} \\
\frac{5}{12}& 0 & \frac{5}{12} &\frac{1}{3}  & \frac{11}{24} & \frac{11}{24} &\frac{17}{24}  & \frac{11}{24} & \frac{17}{24}  &\frac{17}{24}  & \frac{11}{24} & \frac{17}{24}\\
\frac{1}{3} &  \frac{5}{12}  & 0   &\frac{5}{12}  &\frac{17}{24}  & \frac{11}{24} &\frac{17}{24}  & \frac{11}{24} & \frac{17}{24}  & \frac{11}{24} &\frac{17}{24}  & \frac{11}{24}\\
\frac{5}{12} & \frac{1}{3}  & \frac{5}{12} & 0 &\frac{11}{24} &\frac{17}{24}  & \frac{11}{24} & \frac{17}{24}& \frac{11}{24} &\frac{17}{24}  & \frac{11}{24} & \frac{17}{24}\\
\frac{17}{24}  & \frac{11}{24} &\frac{17}{24}  & \frac{11}{24} & 0 &1  & \frac{3}{4} &1 &\frac{1}{2} & 1 & \frac{3}{4} &1\\
 \frac{11}{24} &\frac{17}{24}  & \frac{11}{24} & \frac{17}{24}  & 1 & 0 &  1  &\frac{3}{4} & 1 &\frac{1}{2} & 1 &   \frac{3}{4}\\
\frac{17}{24}  & \frac{11}{24} & \frac{17}{24}  & \frac{11}{24}  &\frac{3}{4}     &1   & 0& 1&\frac{3}{4} &1 &\frac{1}{2} &1\\
 \frac{11}{24} &\frac{17}{24}  & \frac{11}{24} & \frac{17}{24}  &1  & \frac{3}{4}     &1   & 0   &1 & \frac{3}{4}        & 1 & \frac{1}{2}\\
\frac{17}{24}  & \frac{11}{24} & \frac{17}{24}  & \frac{11}{24}     &\frac{1}{2} & 1 & \frac{3}{4} &1&0  & 1 & \frac{3}{4} &1 \\
 \frac{11}{24} &\frac{17}{24}  & \frac{11}{24} & \frac{17}{24}  & 1 &\frac{1}{2} & 1 &   \frac{3}{4}  &1   & 0   &1 & \frac{3}{4}  \\
\frac{17}{24}  & \frac{11}{24} & \frac{17}{24}  & \frac{11}{24}    &\frac{3}{4} &1 &\frac{1}{2}&1 & \frac{3}{4}  &1   & 0   &1 \\
 \frac{11}{24} &\frac{17}{24}  & \frac{11}{24} & \frac{17}{24}  & 1 &   \frac{3}{4}  &1   & \frac{1}{2}    &1  & \frac{3}{4}  &1   & 0\\
\end{array}%
\right],$$
 where $r_{ij}$ denotes resistance distance of two
vertices between $i$ and $j.$
 \hfill$\blacksquare$

\section*{Acknowledgments}
Partially supported by National Natural Science Foundation of
China (Nos. 11401004 and 11471016); Natural Science Foundation of
Anhui Province of China (No. KJ2013B105), Anhui Provincial Natural
Science Foundation (No. 1408085QA03).


\begin{thebibliography}{48}
\expandafter\ifx\csname
natexlab\endcsname\relax\fi
\expandafter\ifx\csname bibnamefont\endcsname\relax
  \fi
\expandafter\ifx\csname bibfnamefont\endcsname\relax
  \fi
\expandafter\ifx\csname citenamefont\endcsname\relax
  \fi
\expandafter\ifx\csname url\endcsname\relax
  \def\url#1{\texttt{#1}}\fi
\expandafter\ifx\csname
urlprefix\endcsname\relax\fi

\providecommand{\eprint}[2][]{\url{#2}}

\bibitem{Bondy1976}
\textcolor{blue}{J. A. Bondy, U. S. R. Murty,  Graph Theory,
Springer, New York, 2008.}

\bibitem{Klein1993}
\textcolor{blue}{D. J. Klein, M. Randi$\acute{c}$, Resistance
distance, J. Math. Chem. 12 (1993) 81-95.}



\bibitem{Bu2014}
\textcolor{blue}{C. Bu, B. Yan, X. Zhou, J. Zhou, Resistance
distance in subdivision-vertex join and subdivision-edge join of
graphs, Linear Algebra Appl. 458 (2014) 454-462.}

\bibitem{Zhou2014}
\textcolor{blue}{ J. Zhou, L. Sun,  W. Wang, C. Bu,  Line star
sets for Laplacian eigenvalues, Linear Algebra Appl. 440 (2014)
164-176.}


\bibitem{Y2013}
\textcolor{blue}{Y. J. Yang, D. J. Klein, A recursion formula for
resistance distances and its applications, Discrete Appl. Math.
161 (2013) 2702-2715.}




\bibitem{Feng2014}
\textcolor{blue}{L. H. Feng, G. Yu, K. Xu, Z. Jiang, A note on the
Kirchhoff index of bicyclic graphs, Ars Comb. 114 (2014) 33-40.}






\bibitem{Ya2014}
\textcolor{blue}{Y. J. Yang, D. J. Klein, Comparison theorems on
resistance distances and Kirchhoff indices of S,T-isomers,
Discrete Applied Mathematics, 175 (2014) 87-93.}



\bibitem{YangJ2014}
\textcolor{blue}{Y. J. Yang, The Kirchhoff index of subdivisions
of graphs, Discrete Appl. Math. 171 (2014)  153-157.}


\bibitem{Wang2012}
 \textcolor{blue}{S. L. Wang, B. Zhou, The signless Laplacian spectra of the
corona and edge corona of two graphs, Linear and Multilinear
Algebra, 61 (2013) 197-204.}

\bibitem{Liu2013}
\textcolor{blue}{X. Liu, P. Lu, Spectra of subdivision-vertex and
subdivision-edge neighbourhood coronae, Linear Algebra Appl. 438
(2013) 3547-3559.}


\bibitem{Lu2013}
\textcolor{blue}{P. Lu,  Y. Miao, Spectra of the
subdivision-vertex and subdivision-edge coronae, arXiv:1302.0457.}

\bibitem{McLeman2011}
\textcolor{blue}{ C. McLeman, E. McNicholas, Spectra of coronae,
Linear Algebra Appl. 435 (2011) 998-1007.}

\bibitem{Gopalapillai2011}
\textcolor{blue}{I. Gopalapillai,  The spectrum of neighborhood
corona of graphs. Kragujevac J. Math. 35 (2011) 493-500.}


\bibitem{Horn1991}
 \textcolor{blue}{R. A. Horn, C. R. Johnson, Topics in Matrix Analysis, Cambridge
University Press, 1991.}


\bibitem{Zhang2013}
 \textcolor{blue}{Z. Zhang, Some physical and chemical indices of clique-inserted
lattices, Journal of Statistical Mechanics: Theory and Experiment,
10 (2013), P10004.}

\bibitem{Liu2014}
 \textcolor{blue}{J. B. Liu,  X. F. Pan,  J. Cao, F. F. Hu,  A note on
¡®some physical and chemical indices of clique-inserted
lattices¡¯, Journal of Statistical Mechanics: Theory and
Experiment, 6 (2014), P06006.}

\bibitem{Liu2015}
\textcolor{blue}{ J. B. Liu, X. F. Pan, Asymptotic incidence
energy of lattices, Physica A: Statistical Mechanics and its
Applications, 422 (2015) 193-202.}

\bibitem{LiuP2015}
\textcolor{blue}{ J. B. Liu, X. F. Pan, F. T. Hu, F. F. Hu,
Asymptotic Laplacian-energy-like invariant of lattices, Appl.
Math. Comput. 253 (2015) 205-214.}



\end{thebibliography}
\end{document}